\documentclass[a4paper,10pt]{article}

\usepackage{amsfonts,amssymb,mathrsfs,amscd}
\def \R {{\mathbb {R}}}
\def \N {{\mathbf {N}}}

\def\eps{\varepsilon}

\begin{document}
\title{\bf Universal Weighted Averages for Ergodic Flows }
\author{\bf Valery V. Ryzhikov }
\date{}

\maketitle
\Large
\begin{abstract}
This paper studies homothetic and more general weighted averages
for flows. Absolutely continuous convolutions of singular weights are considered, thereby strengthening the Kozlov-Treschev result on nonuniform averages for ergodic flows.
The concept of almost mixing, formulated in terms of homothetic averaging, is proposed. An example of a non-mixing flow with almost mixing is given. It is proved that rigid flows do not have this property.

\vspace{2mm}

\end{abstract}

\section{Homothetic and Universal Averages, Weight Convolutions}
This paper considers ergodic flows $T_t$ on the standard probability space $(X,\mu)$ that preserve the measure $\mu$.
In \cite{KT}, homothetic averages of the form
$$P_tf(x)=\int_R f(T_{rt}x) h(r) dr, \ \int_R h\,dr=1, \ h\geq 0.$$
In this case, the convergence of the averages to a constant in the $L_1$ norm was proven, which is a generalization of the classical theorem for ergodic flows, which considers the function $h=\chi_{[0,1]}$.

\vspace{3mm}
\bf Theorem 1.1 (\cite{KT}). \it Let $\nu$ be a normalized measure on $\R$ that is absolutely continuous with respect to Lebesgue measure. Then, for every $f\in L_1(X,\mu)$, $\mu(X)=1$, we have
$$\left\|\int_R f(T_{rt}x) d\nu(r) -\int f \, d\mu\right\|_1\to 0, \ \ t\to +\infty.$$

\rm

\vspace{3mm}
Note that homothetic averages were considered in [cite{91}] for flows with multidimensional time in order to study the multiple mixing property.
In particular, in \cite{91} for an ergodic
flow $T_t$ and measurable $A,B\subset X$ and $h={\chi_D}/{m(D)} $, where
$D \subset \R$ has finite measure, the convergence of
$$\int_\R \mu(A \cap T_{rt}B)h(r) dr\to \mu(A)\mu(B)$$ was shown as $t\to\infty$.

Theorem 1.1 and other results of \cite{KT} stimulated a number of studies on nonuniform averagings (see \cite{B} and references). We call a normalized measure $\nu$ on $R$ universal if the assertion of Theorem 1.1 holds for every ergodic flow and every $ f\in L_1(\mu)$.
The following question arises.

\vspace{3mm}
\it Which singular measures are universal? \rm

\vspace{3mm}
In this note, we show that the class of universal measures includes some singular measures. We will exploit the fact that
universality is preserved under extraction of convolution roots.
Convolution roots of an absolutely continuous measure, which is universal by Theorem 1.1, can be singular. Examples of measures on $\R$ that are singular together with their convolution powers of order $d< n$, but absolutely continuous convolution powers of order $d\geq n$, $n>2$, provide spectra of Sidon automorphisms from \cite{26}. For $n=2$, such measures have long been well known.

Spectral measures of automorphisms are known to be considered on the unit circle in the complex plane. By unfolding the circle into a segment on $\R$, we obtain the desired measure on a line. Alternatively, we can proceed. Using Sidon automorphisms of $S$, we consider a special flow $T_t$ over $S$ with constant return function equal to 1. Alternatively, we construct suitable Sidon flows by analogy with automorphisms. The spectral measures of such flows will be the desired measures on $\R$.

\vspace{3mm}
\bf Theorem 1.2. \it If some convolution power $\nu^{\ast n}$ of a normalized Borel measure $\nu$ on $\R$ is universal for  ergodic flows, then the measure $\nu$ is also universal. \rm

\vspace{3mm}
Open question: \it Does there exist a universal measure for which all convolution powers are singular? \rm

\section{Universal sequence of probability  measures on $\R$}
A sequence of  probability atomless Borel measures $\nu_j$ on $\R$ is called universal if
for every ergodic flow $T_t$ and every function $f\in L_1(\mu)$, we have
$$\left\|P_jf -\int f \, d\mu\right\|_1\to 0, \ \ j\to +\infty,$$
where
$$P_jf(x)=\int_R f(T_{t}x) d\nu_j(t).$$

If the operator $P$ corresponds to a weight $\nu$, then the degree
$P^n$ corresponds to the convolution degree $\nu^{\ast n}=\nu\ast \nu\ast \dots\ast \nu$ ($n$ factors).

Taking the above into account, Theorem 1.2 follows from the following fact.

\vspace{3mm}
\bf Theorem 2.1. \it  Let $\nu_j$ be a sequence of probability Borel measures on $\R$. If for some $n$ the convolution powers of $\nu_j^{\ast n}$ form a universal sequence, then $\nu_j$ is also universal. \rm

\vspace{3mm}
Recall that an operator $Q$ acting on a Hilbert space is called normal if $Q^\ast Q=QQ^\ast$.

\vspace{3mm}
\bf Lemma 2.2. \it Let $Q_j$ be a sequence of normal operators in $L_2$ and for every function $f\in L_2$ with zero average, $\|Q_j^n f\|_2\to 0$. Then $\|Q_j f\|_2\to 0$.\rm

\vspace{3mm}
Proof. A normal operator is equivalent to the operator of multiplication by a bounded function in a suitable measure space. For such operators, the assertion of the lemma is obvious. Since one can do without the spectral representation of normal operators, we present below an elementary proof.

If $\|Q_j^n f\|_2\to 0$, then, since the operators
$Q_j^\ast$ and $Q_j $ commute, we have $$(Q_j^\ast Q_j)^{2^n}f\to_w 0,$$
therefore,
$$\left( Q_j^{2^{n-1}}f\, ,\, Q_j^{2^{n-1}}f\right)\ \to 0.$$
We continue the descent in powers:
$$ (Q_j^\ast Q_j)^{2^{n-2}}\to_w 0,\ \dots, \ (Q_j^\ast Q_j)^2\to_w 0, \ \ Q_j^\ast Q_jf\to_w 0.$$
The last of the convergences is equivalent $\|Q_j f\|_2\to 0$.

Lemma 2.2 can be proved using the spectral representation of a unitary flow. The flow $T_t$ on the space $L_2(\mu)$ is a unitary flow $T_t$; we restrict it to a cyclic space with a cyclic vector $f$, $\|f\|_2=1$. This restriction is isomorphic to the operator of multiplication by $e^{irt}$ in $L_2(\R,\sigma)$ for some probability measure $\sigma$
(called the spectral measure of the flow). The cyclic vector $f$ then corresponds to a constant function ${\bf 1}\in L_2(\R,\sigma)$.
The operator $P_j$ in the spectral representation will correspond to the operator of multiplication by the function
$F(r)=\int_R e^{irt}d\nu_j(t)$ (the conjugate Fourier transform of the measure
$\nu_j$). If $\|F_j(r)^n{\bf 1}\|_2\to 0$ holds for some $n$, then $\|F_j(r){\bf 1}\|_2\to 0$ automatically holds.

\vspace{3mm}
\bf Remark. \rm The recent paper \cite{P} is devoted to the study of the convergence rates of nonuniform averages, where the Fourier transform of probability measures on $\R$ is significantly used.

Proof of Theorem 2.1. The operators $P_j$ act identically on constants. Let $Q_j$ denote the restriction of $P_j$ to the space $L^0_1$ of integrable functions with zero mean.
Bounded functions $g$ with zero mean are dense in $L^0_1$ and in $L^0_2$.
Since
$\|Q_j g\|_1\leq \|Q_j g\|_2\to 0,$ for every $\eps>0$, for $\|f-g\|_1<\eps$, for all sufficiently large $j$,
$\|Q_j f\|_1\leq \eps$ holds. The theorem is proved.

\section{Homothetic Averaging and Almost Mixing}
Recall that a flow $T_u$, $u\in \R$, with an invariant probability measure $\mu$ is mixing if for every $f,g\in L_2(\mu)$ 
$$(T_uf,g)\to \int f \,d\mu \ \int f \,d\mu, \ \ u\to \infty. $$
We call an ergodic flow \it almost mixing \rm if the assertion of Theorem 1.1 holds for any continuous probability Borel measure $\nu$.
The following fact can be considered as a special continuous analog of the well-known theorem from \cite{BH}.

\vspace{3mm}
\bf Theorem 3.1. \it If a flow is mixing,
then it is almost mixing. \rm

\vspace{3mm}
Proof. Let us show that for a continuous measure $\nu$ and a
function $f$ with zero average, $P_t^\ast P_t f\to_w 0$, $t\to+\infty$.
where the operators $P_t$ in $L_2(\mu)$ are defined by the equality
$$P_tg(x)= \int_\R g(T_{tr} x)\,d\nu(r).$$
We have
$$(P_tf,P_tg)= \int_X \int_\R f(T_{tr} x)\,d\nu(r)\ \int_\R g(T_{ts} x)\,d\nu(s) d\mu=$$
$$=\int_X \int_R\int_\R f(T_{t(r-s)} x)\,\, g(x)\,d\nu(r)\, d\nu(s) d\mu.\eqno (1)$$
For $t\to +\infty$
due to the continuity of the measure $\nu$ for a fixed $N$, we have
$$\nu\times\nu( (r,s)\, : \, t(r-s)|<N )\to 0.$$
Since the flow is mixing and $\int f dm=0$, for $\eps>0$ there exists $N$ such that for $u>N$
$$\left| \int_X f(T_u x)\,\, g(x) d\mu \right| <\eps.$$
Thus, we arrive at $|(P_tf,P_tg)|< \eps$ for all sufficiently large $t$.

We have shown that $(P_tf,P_tf)\to 0$, i.e., $\| P_t f\|_2 \to 0.$
The theorem is proved.

We will show that the almost mixing property is not equivalent to the mixing property.

\vspace{3mm}
\bf Theorem 3.2. \it A non-mixing flow can be almost mixing. \rm

\vspace{3mm}
Proof. We use the following assertion (Theorem 2, \cite{21}).
\it If the complement of a set $M\subset \N$ contains segments $[h_j-L_j, a_j+L_j]$ and $h_j,L_j\to \infty$, then the set $M$ is a mixing set for some non-mixing automorphism $S$. \rm

Below, we consider the case where the sequence $h_j$ grows very rapidly, while the sequence $L_j$ grows very slowly.

A flow $T_t$ with similar properties can be obtained as follows. First, a rank-one flow $S_t$ with infinite invariant measure is constructed.
For example, this is a special flow with constant return function over a rank-one automorphism with parameters $r_j=2$, $\bar s_j=(0, s_j(2))$, $s_j(2)/h_j\to \infty$. Such automorphism constructions were used in \cite{21}. The flow $S_t$ is non-mixing: for all sets of finite measure, $\mu(A \cap S_{h_j}A) \to \mu(A)/2.$
However, its main feature is that for every non-mixing sequence $t_i$, there exists a sequence $j_i$ such that $|t_i - h_{j_i}|/h_{j_i} \to 0.$
Next, we consider a Poisson or Gaussian supension $T_t$ over the flow $S_t$, which inherits the mixing properties of the flow $S_t$.

Thus, for $T_t$, deviations from mixing are asymptotically observed only on the intervals $I_j=[h_j-L_j, h_j+L_j]$. These intervals are relatively small and so far apart that, informally speaking, homothetically stretched measures in the limit do not integrally experience the deviation from mixing.

Let a measure $\nu$ be supported in $[0,N]$. Fix a small $\delta>0$.
Consider the intervals $t^{-1}I_j$ for a large $t$.
The interval $[\delta,N]$ cannot contain two intervals of the form $t^{-1}I_j$.
The continuity of the measures defined on the compact set $[0,N]$ (and similarly for $[-N,N]$) implies that the $\nu$-measure of a small interval is small. Therefore, with integration similar to that used in the proof of the
previous theorem (see (1)), deviations from mixing are observed on a set of such small $\nu$-measure that the proof of almost mixing
for a mixing flow also applies to our flow. Finally, the general case easily follows from the above, since $\nu (\R\setminus [-N,N])\to 0, \ N\to\infty$), this completes the proof.

\section{Rigid flows are not almost mixing}
The contrast to mixing is the rigidity of a flow, or the presence of a rigid factor. The latter property averages that
there exists a set $A$, $0<\mu(A)<1$, and a sequence $t(i)\to\infty$ such that
$$\mu(A \cap T_{t(i)}A) \to \mu(A).$$

\vspace{3mm}
\bf Theorem 4.1. \it Every flow with a rigid factor is not almost mixing. \rm

\vspace{3mm}
The theorem implies, for example, that ergodic windings of a torus are not almost mixing. Flows without rigid factors may not be almost mixing. Let us explain. A flow $T_t$ is, by definition, $a$-partially mixing, $a\in (0,1)$, if for every measurable $A,B$, $$\liminf_{t\to\infty} \mu(A \cap T_{t}B) \geq a \mu(A)\mu(B).$$ For any  $a<1$  it is possible to construct examples of partially mixing flows that are not almost mixing.

If significant deviations from mixing are observed on long progressions $t_i, 2t_i, \dots, m(i)t_i$, where $m(i)\to\infty$, $t_i\to\infty$, then 
almost-mixing fails.

Proof of the theorem.
From the convergence of $\mu(A\,\Delta\, T_{t(i)}A) \to 0$ it obviously follows that
$$ \mu(A\,\Delta\, T_{2t(i)}A) \to 0,\ \ \mu(A\,\Delta \,T_{3t(i)}A) \to 0, \ \dots, \
\mu(A\,\Delta\, T_{m(i)t(i)}A) \to 0$$
for some sequence $m(i)\to+\infty$. \rm

Fixing the specified set $A$, we define a measure $\nu$ such that for
homothetic averages
$$P_{s}f=\int_R f(T_{sr}x) d\nu(r)$$
with $$s(i):=m(i)t(i)\to\infty$$
we have
$$(P_{s(i)} \chi_A, \chi_A) \to \mu(A)\neq \mu(A)^2.$$
This will demonstrate the absence of almost mixing for the flow $T_t$.

The measure $\nu$ is defined as follows: its support lies in $[0,1]$, in two disjoint segments, and in each of these segments, the support again lies in two disjoint smaller segments. We call such segments level-2 segments. In each of these four segments, the support of the measure is again located in two very small level-3 segments, and so on.
The measure $\nu$ is such that each of the $2^n$ segments of level $n$ has $\nu$-measure equal to $2^{-n}$.
The set of points that is the intersection of the unions of segments of the same level, similar to the standard Cantor set. The measure $\nu$, by construction,
is uniformly distributed on this set.

Now we indicate how these segments are chosen. Consider one of the segments of level $n-1$, denoting it by $[a,b]$.
We choose disjoint segments $I, I'\subset [a,b]\subset [0,1]$
as follows. For all sufficiently large $i$, the segments
$[s(i)a,\, s(i)b]$ contain the points $p t(i)$ and $(p+1)t_i$, where $p+1<m(i)$.

Let $\delta_n >0$ be such that
$$|(\chi_{T_tA}, \chi_A) - \mu(A)|<\frac 1 n$$
for
$$t\in [p t(i)-\delta_n, \, p t(i)+\delta_n]\ \cup \ [(p+1) t(i)-\delta_n, \,
(p+1) t(i)+\delta_n].$$
Set (recall that $s(i)=m(i)t(i)$)
$$I=\left[\frac {p t(i)-\delta_n}{s(i)}\, , \, \frac {p t(i)+\delta_n}
{s(i)}\right], \ \
I'=\left[\frac {(p+1) t(i)-\delta_n}{s(i)}, \,\frac {(p+1) t(i)+\delta_n}{s(i)}\right].$$

Thus, we determine the required intervals $I,I'$ of level $n$ in each of the intervals of level $n-1$. Moreover, for all intervals of level $n$, a common $i$ can be easily selected,  we denote it by $i_n$.

For $\nu$-almost all $r\in\R$ one has $$\mu(T_{s(i_n)r}A\, \Delta\, A)\to 0, \ \ n\to\infty,$$
therefore
$$\int_X \int_R \chi_A(T_{s(i_n)r}x)\, \chi_A \, d\nu(r)\,d\mu =
(P_{s(i_n)}\chi_A\, ,\, \chi_A)\ \to \mu(A), \ \ n\to\infty,$$
what completes the proof.

\end{document}